\documentclass[12pt,psamsfonts]{amsart}

\usepackage{amsmath,amssymb}
\usepackage{graphicx}
\usepackage{color}
\usepackage{psfrag} 
\usepackage{overpic} 

\newtheorem{theorem}{Theorem}
\newtheorem{proposition}[theorem]{Proposition}

\newtheorem{remark}[theorem]{Remark}
\newtheorem{definition}[theorem]{Definition}

\title[Cyclicity of Rigid centers on center manifolds]{Cyclicity of Rigid centers on center manifolds of three-dimensional systems}
\author[C. Pessoa, L. Queiroz, J. Ribeiro]{}


\begin{document}
 \maketitle

\centerline{\scshape  Claudio Pessoa,  \; Lucas Queiroz, \; Jarne Ribeiro}
\medskip

{\footnotesize \centerline{Universidade Estadual Paulista (UNESP), Instituto de Bioci\^encias Letras e Ci\^encias Exatas,} \centerline{R. Cristov\~ao Colombo, 2265, 15.054-000, S. J. Rio Preto, SP, Brasil }
\centerline{\email{c.pessoa@unesp.br} and \email{lucas.queiroz@unesp.br}}}

{\footnotesize \centerline{Instituto Federal de Educa\c{c}\~ao, Ci\^encia e Tecnologia do Sul de Minas Gerais - IFSULDEMINAS,} \centerline{R. Mario Ribola, 409, Penha II, 37903-358, Passos, MG, Brasil }
	\centerline{\email{jarne.ribeiro@ifsuldeminas.edu.br}}}

\medskip

\bigskip

\begin{quote}{\normalfont\fontsize{8}{10}\selectfont
{\bfseries Abstract.}
We work with polynomial three-dimensional rigid differential systems. Using the Lyapunov constants, we obtain lower bounds for the cyclicity of the known rigid centers on their center manifolds. Moreover, we obtain an example of a quadratic rigid center from which is possible to bifurcate 13 limit cycles, which is a new lower bound for three-dimensional quadratic systems.
\par}
\end{quote}

\section{Introduction}

In 1900, the International Mathematical Congress took place in Paris. In this event, the acclaimed mathematician David Hilbert proposed a list of open problems among which is the problem of determining and locating the limit cycles for a given planar polynomial differential system. This is known as \emph{Hilbert's 16th problem} and it is considered to be the hardest challenge within the Qualitative Theory of Differential Equations, since it remains open even for the simpler cases, i.e. quadratic polynomial vector fields (see \cite{ChavarrigaSurvey, JibinLi}). There are several weaker versions of this problem, for instance, determining the quantity $M(n)$ of limit cycles that bifurcate from an elementary center or focus for polynomial vector fields of degree $n$ \cite{ChristopherLi, Romanovski, Zoladek1}. A singular point of a planar polynomial vector field is an \emph{elementary center or focus} if the eigenvalues of its Jacobian matrix evaluated at the singular point are purely imaginary. Any differential system associated to a vector field having an elementary center or focus can be written, after the proper change of variables and time rescaling, in the form:
\begin{equation*}\label{eqPlanar1}
\begin{array}{lcr}
\dot{x}=-y+P(x,y),\\
\dot{y}=x+Q(x,y),
\end{array}
\end{equation*}
where $P$ and $Q$ are polynomials with no linear or constant terms.

The standard method to estimate the quantity of limit cycles bifurcating from an elementary center or focus is to compute the Lyapunov coefficients (or equivalently the focal coefficients) when small perturbations of the system are considered (see \cite{Bautin,Christopher,Romanovski} and the references therein). Christopher \cite{Christopher} developed a simple approach to obtain lower bounds for the quantity of limit cycles bifurcating from an elementary center or focus, i.e. the small amplitude limit cycles. The idea consists in considering the linear part of the power series expansion of the first $k$ Lyapunov coefficients. If those $k$ linear parts are linearly independent, then it is possible by small perturbations to obtain $k$ limit cycles, that is to say that the cyclicity of the considered singular point is at least $k$. 

We consider three-dimensional differential systems with a Hopf singular point, that is, systems having a singular point for which the Jacobian matrix has a pair of purely imaginary eigenvalues and a non-zero real eigenvalue. These systems can be put in the following canonical form, by means of a linear change of variables and time rescaling:
\begin{equation}\label{eq1}
	\begin{array}{lcr}
		\dot{x}=-y+P(x,y,z),\\
		\dot{y}=x+Q(x,y,z),\\
		\dot{z}=-\lambda z+R(x,y,z),
	\end{array}
\end{equation}
where $P,Q,R$ are polynomials with no linear nor constant terms and $\lambda\neq 0$. The Center Problem and Hilbert's 16th can be naturally stated for vector fields \eqref{eq1}, since there exists an invariant bidimensional $C^r$-manifold $W^c$ tangent to the $xy$-plane at the origin for every $r\geqslant 3$. This result is known as the \emph{Center Manifold Theorem} (see \cite{Kelley, Sijbrand}). 

Recently, these problems have attracted the attention of researchers in the qualitative theory of ordinary differential equations. In \cite{Adam}, the authors prove that the set of systems of form \eqref{eq1} having a center on the local center manifold at the origin, to a fixed value of $\lambda$, corresponds to a variety in the space of admissible coefficients. In \cite{mahdi2011centers} the authors solve the Center Problem for the Lü system and Mahdi in \cite{MahdiJerk} considers the Center Problem for quadratic systems of the form \eqref{eq1} obtained from a third order ODE. The authors of \cite{MahdiPessoa} complete this study using a new hybrid symbolic-numerical approach.

Regarding the cyclicity problem, the results of Christopher are originally stated for planar systems, but they can be extended for system \eqref{eq1} (see \cite{Garcia}), since, for these systems, the Lyapunov coefficients are well-defined and play the same role as their bidimensional counterpart \cite{Adam, Garcia}. In \cite{garcia2019center,QueirozGouveia,Ivan,SangFercec,PeiYu}, the authors use these ideas to study the cyclicity of centers on center manifolds of several families of three-dimensional systems.

In the paper \cite{mahdi_pessoa_ribeiro_2021}, the authors extended the concept of rigid systems and rigid centers for three-dimensional systems \eqref{eq1}. In essence, system \eqref{eq1} has a rigid center at the origin when its restriction to the center manifold has a rigid center at the origin. Moreover, in \cite{mahdi_pessoa_ribeiro_2021}, the authors solved the Center Problem for several families of rigid systems.

Our goal in this work is to obtain lower bounds for the quantity of small amplitude limit cycles which bifurcate from the rigid centers identified in \cite{mahdi_pessoa_ribeiro_2021}. The structure of this paper is as follows: In Section 2, we present the preliminary results on the investigation of the cyclicity for three-dimensional systems as well as an overview of the concept of rigid systems studied in \cite{mahdi_pessoa_ribeiro_2021}. In Section 3 we study the cyclicity of the rigid centers considered in \cite{mahdi_pessoa_ribeiro_2021} obtaining lower bounds for the cyclicity of each system under quadratic perturbations. Also in this section, we prove the main result of this work in which we present an example of a quadratic three-dimensional system from which is possible to obtain $13$ small limit cycles from quadratic perturbations. This is the highest known lower bound for quadratic three-dimensional systems.

\section{Fundamental Concepts}

\subsection{Lyapunov coefficients and limit cycles}

We work with the following perturbation of system \eqref{eq1}:
\begin{equation}\label{eqPerturbation}
\begin{array}{lcr}
		\dot{x}=-y+P(x,y,z)+G_1(x,y,z),\\
		\dot{y}=x+Q(x,y,z)+G_2(x,y,z),\\
		\dot{z}=-\lambda z+R(x,y,z)+G_3(x,y,z),
	\end{array}
\end{equation}
where $P,Q,R,G_1,G_2,G_3$ are polynomials with no constant nor linear terms. Let $\Lambda=(\lambda_1,\dots,\lambda_m)$ denote the perturbation parameters, i.e. the coefficients of $G_i$, $i=1,2,3$. We assume that for $\Lambda=0$, the origin is a center on the center manifold of system \eqref{eqPerturbation}.

First we need to define a \emph{displacement map} to search for limit cycles. Via the change of variables $x=\rho\cos\theta$, $y=\rho\sin\theta$ and $z=\rho\omega$, the solution curves of \eqref{eqPerturbation} are determined by the following equations:
\begin{equation}\label{eqRW}
\begin{array}{lcr}
\dfrac{d\rho}{d\theta}=\dfrac{\psi(\rho,\theta,\Lambda)}{1+\phi(\rho,\theta,\Lambda)},\\
\dfrac{d\omega}{d\theta}=\dfrac{-\mu\omega+\Omega(\rho,\theta,\Lambda)}{1+\phi(\rho,\theta,\Lambda)}.
\end{array}
\end{equation}
For each $(\rho_0,\omega_0)$ with sufficiently small $\|(\rho_0,\omega_0)\|$, let $\varphi(\theta,\rho_0,\omega_0)=(\rho(\theta,\rho_0,\omega_0),\omega(\theta,\rho_0,\omega_0))$ be the solution of \eqref{eqRW} with initial conditions $(\rho(0,\rho_0,\omega_0),\omega(0,\rho_0,\omega_0))=(\rho_0,\omega_0)$.
\begin{definition}For system \eqref{eq1}, the map $$d(\rho_0,\omega_0)=\varphi(2\pi,\rho_0,\omega_0)-(\rho_0,\omega_0)$$ is called \emph{displacement map}.\end{definition}

The displacement map $d(\rho_0,\omega_0)=(d_1(\rho_0,\omega_0),d_2(\rho_0,\omega_0))$ is a bidimensional map, unlike the planar counterpart. However, there exists a unique analytical function $\tilde{\omega}(\rho_0)$ defined in a neighborhood $V$ of $\rho_0=0$ such that $d_2(\rho_0,\tilde{\omega}(\rho_0))\equiv 0$ (see \cite{GarciaJacobi}). The analytical function $\mathbf{d}(\rho_0)=d_1(\rho_0,\tilde{\omega}(\rho_0))$, is called  \emph{reduced displacement map}. We write its power series as:
$$\mathbf{d}(\rho_0)=v_1(\Lambda)\rho_0+v_2(\Lambda)\rho_0^2+v_3(\Lambda)\rho_0^3+v_4(\Lambda)\rho_0^4+\cdots.$$
It is known that the first nonzero coefficient $v_k(\Lambda)$ is the coefficient of an odd power of $\rho_0$ (see \cite{GarciaJacobi}). Thus, we define the coefficients $l_k=v_{2k+1}(\Lambda)$ as the \emph{Lyapunov coefficients}. The isolated zeros of the reduced displacement map correspond to limit cycles of system \eqref{eqPerturbation}. Furthermore, the origin is a center on the center manifold $W^c$ if and only if $l_k=0$ for all $k\in\mathbb{N}$. 

For computational purposes, working with the Lyapunov coefficients is not very efficient. Thus, another method is applied to study the displacement function.

Following the steps of \cite{Queiroz, Adam, Garcia}, we construct a formal series $$H(x,y,z)=x^2+y^2+\sum_{j+k+l\geqslant 3}p_{jkl}x^jy^kz^l,$$
with real coefficients $p_{jkl}$ yet to determine. Let $X$ denote the vector field associated to system \eqref{eqPerturbation}. It is always possible to choose $p_{jkl}$ such that
\begin{equation}\label{XH=Lyap}
XH=\sum_{k\geqslant2}L_{k-1}(x^2+y^2)^k.
\end{equation}
Moreover, the quantities $L_{k}$ are rational functions on the parameters of system \eqref{eqPerturbation} including $\Lambda$. The coefficients $L_{k}$ are called \emph{focal coefficients} for system \eqref{eq1}. The focal coefficients and the Lyapunov coefficients have the following relationship: given a positive integer $k>2$, we have
$$l_1=\pi L_{1}\;\mbox{ and }\;l_{k}=\pi L_{k} \mbox{\hspace{0.3cm} \emph{mod}}\langle L_{1},L_{2},\dots,L_{k-1}\rangle.$$
This result is proven in \cite{Garcia}. Since the computations involving the expressions given in \eqref{XH=Lyap} are algebraic, one can use symbolic mathematics softwares (e.g. Maple and Mathematica) to compute a large amount of focal coefficients in less computational time.


Once a sufficient amount of focal coefficients is computed, we can use the next result to estimate the cyclicity of the center of \eqref{eqPerturbation}$\mid_{\Lambda=0}$.

\begin{theorem}\label{TeoBifLinear}
For system \eqref{eqPerturbation}, if
$${\rm rank}\left[\dfrac{\partial (L_{1},L_{2},\dots,L_{k})}{\partial (\lambda_1,\dots,\lambda_m)}\right]_{\Lambda=0}=k,$$
then it is possible to choose parameters $\Lambda$ such that at least $k-1$ small amplitude limit cycles bifurcate from the origin. 
\end{theorem}

The proof of this theorem can be found in \cite{Garcia}. It is a simple computational task to compute the matrix in Theorem \ref{TeoBifLinear}. In fact, since the perturbation of the center in system \ref{eqPerturbation} is given through the functions $G_i$, $i=1,2,3$, it is sufficient to compute the linear term of the Taylor expansion at $\Lambda=0$ of each focal coefficient. The idea of using the linear parts of the Lyapunov coefficients (or focal coefficients) to estimate the cyclicity of centers was first introduced by Christopher in \cite{Christopher} for the planar case.

It is known in the literature (see, for instance, \cite{Garcia}) that an additional limit-cycle can be obtained by making a perturbation of the trace of the linear part of system \eqref{eq1}, i.e. the following perturbation:
\begin{equation}\label{eqPerturbationtrace}
	\begin{array}{lcr}
		\dot{x}=\lambda_0x-y+P(x,y,z)+ G_1(x,y,z),\\
		\dot{y}=x+\lambda_0y+Q(x,y,z)+ G_2(x,y,z),\\
		\dot{z}=-\lambda z+R(x,y,z)+ G_3(x,y,z).
	\end{array}
\end{equation}

In \cite{TorreGouveia}, the authors proposed an extension of Christopher's ideas. Their approach is based in the next theorem, proven in \cite{TorreGouveia}, which consists in the study of the higher order developments of the Lyapunov coefficients. 

\begin{theorem}\label{TeoBifHighOrder}
	Suppose that system \eqref{eqPerturbation} satisfies the hypothesis of Theorem \ref{TeoBifLinear}. Then, by a change of variables if necessary, we can assume that $l_i=u_i+O(|u|^2)$ for $i=0,\dots,k$ and the next Lyapunov coefficients $l_i=h_i(u)+O(|u|^{m+1})$, for $i=k+1,\dots,k+l$, where $h_i$ are homogeneous polynomials of degree $m\geqslant 2$ and $u=(u_{k+1},\dots, u_{k+l})$. If there exists a line $\eta$, in the parameter space, such that $h_i(\eta)=0,$ for $i=k+1,\dots,k+l-1$, the hypersurfaces $h_i=0$ intersect transversally along $\eta$ for $i=k+1,\dots, k+l-1$, and $h_{k+l}(\eta)\neq 0$, then there are perturbations of the center which produce $k+l$ limit cycles.
\end{theorem}

Theorem \ref{TeoBifHighOrder} has been successfully applied to obtain new lower bounds for the cyclicity of several polynomial systems, for instance in the papers \cite{QueirozGouveia, TorreGouveia}.

\subsection{Rigid Systems in $\mathbb{R}^3$}

We exhibit a quick overview of rigid systems in $\mathbb{R}^3$, which are studied with more detail in \cite{mahdi_pessoa_ribeiro_2021}.

\begin{definition}[Rigid systems in $\mathbb{R}^3$]
We say that the three-dimensional system \eqref{eq1} is rigid when its restriction to a center manifold is rigid, i.e., the restricted system has constant angular speed. Moreover, if the origin is a center on the center manifold, we say that the origin is a rigid center.
\end{definition}

In order to avoid restricting system \eqref{eq1} to a center manifold, which is not effective computationally, the authors in \cite{mahdi_pessoa_ribeiro_2021} introduced a subclass of rigid systems.

\begin{definition}
We say that system \eqref{eq1} is rigid by cylindrical coordinates when, in cylindrical coordinates $(x,y,z)\to (\rho\cos\theta,\rho\sin\theta,z)$, the transformed system satisfies $\dot{\theta}=1$.
\end{definition}

The orbits of rigid systems by cylindrical coordinates rotate around the $z$-axis with constant angular speed. The following result classify rigid systems \eqref{eq1}.

\begin{proposition}[Corollary 2.5 in \cite{mahdi_pessoa_ribeiro_2021}]
Consider system \eqref{eq1}. Then,
\begin{itemize}
\item[(a)] it is a rigid system if and only if $$\left(xQ(x,y,z)-yP(x,y,z)\right)\vert_{z=h(x,y)}\equiv 0,$$ where $h(x,y)$ is the local expression of a local center manifold;
\item[(b)] it is a rigid system by cylindrical coordinates if and only if $xQ(x,y,z)-yP(x,y,z)\equiv 0$.
\end{itemize}
\end{proposition}

\section{Cyclicity of Rigid Centers in $\mathbb{R}^3$}


Consider the following rigid system:
\begin{equation}\label{eqRigidnm}
	\begin{array}{lcr}
		\dot{x}=-y+xF_n(x,y,z),\\
		\dot{y}=x+yF_n(x,y,z),\\
		\dot{z}=-\lambda z+R_m(x,y),
	\end{array}
\end{equation}
where $F_n$ is a homogeneous polynomial of degree $n$ and $R_m$ is a homogeneous polynomial of degree $m$. The parameter values for which the above system has a rigid center at the origin are called \emph{center conditions} and were obtained in \cite[Theorem 4.1]{mahdi_pessoa_ribeiro_2021}. For each center condition we perform the perturbation \eqref{eqPerturbation}, with $G_i,i=1,2,3$ being quadratic polynomials to study the cyclicity. More precisely, we obtained the following result.

\begin{theorem}\label{TeoRigidnm}
Consider system \eqref{eqRigidnm} with $F_n(x,y,z)=\sum_{j+k+l=n}a_{jkl}x^jy^kz^l$ and $R_m(x,y)=\sum_{j+k=m}b_{jk}x^jy^k$.
\begin{itemize}
\item[a)] For $n=1$ and $m=2$, the two center conditions are $\{a_{001}=0\}$ and $R_2\equiv 0$. The cyclicity of the center is at least $6$ for the first center condition and at least $4$ for the last;
\item[b)] For $n=1$ and $m=3$, the two center conditions are $\{a_{001}=0\}$ and $R_3\equiv 0$. The cyclicity of the center is at least $7$ for the first center condition and at least $4$ for the last;
\item[c)] For $n=1$ and $m=4$, the two center conditions are $\{a_{001}=0\}$ and $R_4\equiv 0$. The cyclicity of the center is at least $8$ for the first center condition and at least $4$ for the last;
\item[d)] For $n=2$ and $m=2$, with the center condition $\{a_{200}=-a_{020}, a_{101}=a_{011}=a_{002}=0\}$, the cyclicity of the center is at least $4$.
%
\end{itemize}
\end{theorem}
\noindent\textbf{Proof: }Consider the following perturbation of system \eqref{eqRigidnm}:
\begin{equation*}
	\begin{array}{lcr}
		\dot{x}=-y+xF_n(x,y,z)+G_1(x,y,z),\\
		\dot{y}=x+yF_n(x,y,z)+G_2(x,y,z),\\
		\dot{z}=-\lambda z+R_m(x,y)+G_3(x,y,z),
	\end{array}
\end{equation*}
with $G_i$, $i=1,2,3$ being homogeneous quadratic polynomials. For each center condition described above, we compute the first 12 focal coefficients of the above perturbed system and the rank of their linear parts with respect to the coefficients of $G_i$, the perturbation parameters. Hence, applying Theorem \ref{TeoBifLinear} and performing a perturbation of the trace, the result holds. For instance, for the first center condition in statement (a), the rank of the linear part of the focal coefficients is $6$. For the other statements, we have: 
\begin{itemize}
\item[a)] Case $n=1,m=2$:
\subitem $\bullet$ $a_{001}=0$: rank $6$;
\subitem $\bullet$ $R_2\equiv 0$: rank $4$.
\item[b)] Case $n=1,m=3$:
\subitem $\bullet$ $a_{001}=0$: rank $7$;
\subitem $\bullet$ $R_3\equiv 0$: rank $4$.
\item[c)] Case $n=1,m=4$:
\subitem $\bullet$ $a_{001}=0$: rank $8$;
\subitem $\bullet$ $R_4\equiv 0$: rank $4$.
\item[d)] Case $n=2,m=2$:
\subitem $\bullet$ $a_{200}=-a_{020}, a_{101}=a_{011}=a_{002}=0$: rank $4$.
\end{itemize}
By Theorem \ref{TeoBifLinear}, the result holds. \qed
\medskip

\begin{remark}
For case $n=2$ and $m=2$ in statement (d) of Theorem \ref{TeoRigidnm}, there is another center condition given by  $\{a_{200}=-a_{020}, b_{20}=b_{11}=b_{02}=0\}$. For this center condition, we were not able to obtain lower bounds for the cyclicity using quadratic perturbations since the linear parts of the focal coefficients are all null. The same situation occurs for the case $n=2$ and $m=3$ with its respective center conditions described in \cite[Theorem 4.1]{mahdi_pessoa_ribeiro_2021}. We believe that it is possible to obtain better lower bounds than those obtained in Theorem \ref{TeoRigidnm} by using higher order perturbations. For instance, for the case $n=2$ and $m=3$ with center condition $\{a_{200}=-a_{020}, a_{101}=a_{011}=a_{002}=0\}$, using cubic perturbations, the cyclicity is at least $15$. Here, we chose quadratic perturbations to simplify the computations since, even with cubic perturbations, the computation of the linear parts of the focal coefficients is computationally a hard task.
\end{remark}

The rigid system
\begin{equation}\label{eqRigidnmcomp}
	\begin{array}{lcr}
		\dot{x}=-y+xF_n(x,y,z),\\
		\dot{y}=x+yF_n(x,y,z),\\
		\dot{z}=-\lambda z+R_m(x,y,z),
	\end{array}
\end{equation}
where $F_n$ is a homogeneous polynomial of degree $n$ and $R_m$ is a homogeneous polynomial of degree $m$ has also been considered in \cite[Theorems 4.3 and 4.4]{mahdi_pessoa_ribeiro_2021} where they obtained center conditions for particular values of $m$ and $n$.
\begin{theorem}\label{TeoCycF1R2}
Consider system \eqref{eqRigidnmcomp} with $F_n(x,y,z)=\sum_{j+k+l=n}a_{jkl}x^jy^kz^l$ and $R_m(x,y,z)=\sum_{j+k+l=m}b_{jkl}x^jy^kz^l$.
For $n=1$ and $m=2$, the known center conditions are:
\begin{itemize} 
\item[a)] $\{a_{001}=0\}$;\medskip
\item[b)] $R_2(x,y,0)\equiv 0$;\medskip
\item[c)] $\{a_{100}=a_{010}=0$, $b_{101}=b_{011}=0$, $b_{002}=2a_{001}$, $b_{020}=-b_{200}\}$;\medskip
\item[d)] $\{b_{101}=2a_{100}+a_{010}$, $b_{020}=-b_{200}$, $b_{011}=-a_{100}+2a_{010}, \linebreak b_{002}=2a_{001}\}$;\medskip
\item[e)] $\{b_{101}=2a_{100}$, $b_{020}=-b_{200}$, $b_{011}=2a_{010}$, $b_{002}=2a_{001},\linebreak (a_{100}^2-a_{010}^2)(b_{200}-b_{110})+a_{100}a_{010}(b_{110}+4b_{200})=0\}$;\medskip
\item[f)]$\left\{b_{200}=\dfrac{a_{100}a_{010}(b_{002}-a_{001})}{a_{001}^2}, b_{110}=\dfrac{(a_{100}^2-a_{010}^2)(a_{001}-b_{002})}{a_{001}^2},\right. \linebreak\left. b_{101}=\dfrac{b_{002}(a_{100} + a_{010})-a_{010}a_{001}}{a_{001}}, b_{020}=\dfrac{a_{100}a_{010}(a_{001}-b_{002})}{a_{001}^2}, \right. \linebreak\left. b_{011}=\dfrac{b_{002}(a_{010}-a_{100}) + a_{100}a_{001}}{a_{001}}\right\}$.
\end{itemize}
The cyclicity of the center at the origin of system \eqref{eqRigidnmcomp} is at least $8$ for condition (a), at least $5$ for condition (b), at least $3$ for condition (c), at least $9$ for condition (d), at least $10$ for condition (e) and at least $11$ for condition (f).
\end{theorem}

\begin{theorem}\label{Teo10}
Consider system \eqref{eqRigidnmcomp} with $F_n(x,y,z)=\sum_{j+k+l=n}a_{jkl}x^jy^kz^l$ and $R_m(x,y,z)=\sum_{j+k+l=m}b_{jkl}x^jy^kz^l$.
For $n=m=2$ and $F_2=R_2$, the two center conditions are $\{a_{200}=a_{110}=a_{020}=0\}$ and $\{a_{020}=-a_{200}, a_{101}=a_{011}=a_{002}=0\}$. The cyclicity of the center is at least $3$ for the first center condition and at least $2$ for the last;
\end{theorem}

\noindent\textbf{Proofs of Theorems \eqref{TeoCycF1R2} and \eqref{Teo10}: }The proofs are analogous to that of Theorem \ref{TeoRigidnm}. More precisely, for each center condition, we compute the first 12 focal coefficients and using Theorem \ref{TeoBifLinear} and performing a perturbation of the trace, we estimate the cyclicity by computing the rank of their linear parts. \qed
\medskip
%
%

The following rigid system
\begin{equation}\label{eqFz}
	\begin{array}{lcr}
		\dot{x}=-y+xF(z),\\
		\dot{y}=x+yF(z),\\
		\dot{z}=-\lambda z+R_2(x,y,z),
	\end{array}
\end{equation}
where $F(z)=\sum_{i=1}^{9}a_iz^i$ and $R_2$ is a homogeneous polynomial of degree $2$ was studied in \cite[Theorem 5.1]{mahdi_pessoa_ribeiro_2021}. 


\begin{theorem}
Consider system \eqref{eqFz} with $F(z)=\sum_{i=1}^{9}a_iz^i$ and $R_2(x,y,z)=\sum_{j+k+l=2}b_{jkl}x^jy^kz^l$. For the three center conditions: $F(z)\equiv 0$, $R_2(x,y,0)\equiv 0$ and $\{a_2=a_4=a_5=\dots=a_9=0, b_{101}=b_{011}=a_3=0$,$b_{002}=2a_1, b_{020}=-b_{200}\}$, the cyclicity of the center is at least $5$, $3$ and $3$, respectively for each center condition.
\end{theorem}
\noindent\textbf{Proof: }Analogous to the proof of the previous theorems. For the first and third center conditions, we compute the first 12 focal coefficients as well as the rank of their linear parts. For the second, we compute the first 10 focal coefficients. The conclusion follows from Theorem \ref{TeoBifLinear} and by performing a perturbation of the trace. \qed

\medskip

Finally, we consider the following rigid system:
\begin{equation}\label{eqFz3R3}
	\begin{array}{lcr}
		\dot{x}=-y+xF(z),\\
		\dot{y}=x+yF(z),\\
		\dot{z}=-\lambda z+R_3(x,z),
	\end{array}
\end{equation}
where $F(z)=a_1z+a_2z^2+a_3z^3$ and $R_3$ is a homogeneous polynomial of degree $3$ which was studied in \cite[Theorem 5.2]{mahdi_pessoa_ribeiro_2021}. 

\begin{theorem}
Consider system \eqref{eqFz3R3} with $R_m(x,z)=\sum_{j+l=3}b_{jl}x^jz^l$. For the two center conditions: $\{b_{30}=0\}$ and $\{a_2=3a_1^2$, $a_3=0$, $b_{21}=b_{12}=0, b_{03}=9a_1^2\}$, the cyclicity of the center is at least $3$ for the first center condition and at least $14$ for the last.
\end{theorem}

\noindent\textbf{Proof: }We proceed analogously to the previous theorems. For the first center condition, we compute the first 12 focal coefficients as well as the rank of their linear part. For the second center condition, we compute 17 focal coefficients. Once again, the conclusion follows from Theorem \ref{TeoBifLinear} and by performing a perturbation of the trace.\qed
\begin{remark}
In \cite[Theorem 5.2]{mahdi_pessoa_ribeiro_2021}, the authors also present the center condition for system \eqref{eqFz3R3} given by $F(z)\equiv 0$. Since under quadratic perturbations, the linear parts of the focal coefficients are all null, we were not able to obtain a lower bound for the cyclicity.
\end{remark}

\medskip


In \cite{Ivan}, the authors obtained a lower bound for the number of small limit cycles bifurcating from the center for polynomial three-dimensio- nal systems with a Hopf singular point for degrees $2$ through $5$. For quadratic systems, their lower bound was $11$ limit cycles. In \cite{QueirozGouveia}, the authors applied the study of the higher order terms of the focal coefficients to find an example of a quadratic system that bifurcates, under quadratic perturbations, 12 limit cycles from the center. Using the same technique for system \eqref{eqRigidnmcomp} with $n=1$ and $m=2$, which consists in applying Theorem \ref{TeoBifHighOrder}, we were able to surpass this lower bound, as we prove the next result.

\begin{theorem}\label{Teo13limitcycles}
Consider the rigid system
\begin{equation}\label{eqRigid12}
	\begin{array}{lcr}
		\dot{x}=-y+xF_1(x,y,z),\\
		\dot{y}=x+yF_1(x,y,z),\\
		\dot{z}=-\lambda z+R_2(x,y,z),
	\end{array}
\end{equation}
with $F_1(x,y,z)=x+2y+3z$ and $R_2(x,y,z)=\tfrac{2}{3}x^2+xy+4xz-\tfrac{2}{3}y^2+3yz+6z^2$. The origin is a rigid center on the center manifold and at least 13 limit cycles bifurcate from it under quadratic perturbations. 
\end{theorem}
\noindent\textbf{Proof: }Note that system \eqref{eqRigid12} satisfies the center condition described in item (d) from Theorem \ref{TeoCycF1R2} and therefore the origin is a center on the center manifold. Considering a quadratic perturbation \eqref{eqPerturbation} with 
$$G_1=\sum_{j+k+l=2}a_{jkl}x^jy^kz^l,\;G_2=\sum_{j+k+l=2}b_{jkl}x^jy^kz^l,\;G_3=\sum_{j+k+l=2}c_{jkl}x^jy^kz^l$$
for system \eqref{eqRigid12}, we compute the first 13 focal coefficients. The rank of their linear part is 9 and by Theorem \ref{TeoBifLinear}, it is possible to obtain 9 limit cycles bifurcating from the origin using the perturbation parameters $(a_{jkl},b_{jkl},c_{jkl})$ as well as a perturbation of the trace \eqref{eqPerturbationtrace}. To simplify the computations, we assign $c_{011}=c_{020}=c_{101}=c_{110}=c_{200}=0$ since these parameters do not play a role when we apply Theorems \ref{TeoBifLinear} and \ref{TeoBifHighOrder} to obtain a higher estimate for the cyclicity. Due to the size of the expressions, we only exhibit the linear part of $L_1,\dots, L_3$:
\begin{eqnarray}
L_1^1&=&\tfrac {22}{45}a_{{011}}+4a_{{020}}-{\tfrac {4}{45}}
a_{{101}}+\tfrac{2}{3}a_{{110}}+\tfrac{4}{3}a_{{200
}}+\tfrac {4}{45}b_{{011}}-\tfrac{2}{3}b_{{020}}\nonumber\\
&&+\tfrac 
{22}{45}b_{{101}}-\tfrac{4}{3}b_{{110}}-2b_{{200}},\nonumber\\
L_2^1&=&\tfrac {3713}{25}b_{{200}}+\tfrac{4}{3}c_{{002}}-\tfrac{2}{9}a_{{002}}-\tfrac {39937}{1125}a_{{011}}-\tfrac {7297}{25}a_{{020}}+\tfrac {6409}{1125}a_{{101}}\nonumber\\
&&-\tfrac {3389}{75}a_{{110}}-\tfrac {1463}{15}a_{{200}}-\tfrac {9409}{1125}b_{{011}}+\tfrac {3839}{75}b_{{020}}-\tfrac {37687}{1125}b_{{101}}+\tfrac {
		1433}{15}b_{{110}},\nonumber\\
L_3^1&=&\tfrac {44084}{1785}a_{{002}}+\tfrac {28460414}{7875}a_{
	{011}}+\tfrac {442099341}{14875}a_{{020}}-\tfrac {76527541}{133875}a_{{101}}-
\tfrac {5168}{35}c_{{002}}\nonumber\\
&&+\tfrac {203415874}{44625}a_{{110}}+\tfrac {444002837}{44625}a_{{200}}+
\tfrac {5707}{1785}b_{{002}}+\tfrac {116062741}{133875}b_{{011}}\nonumber\\
&&-\tfrac {233857774}{44625}b_{{020}}+
	\tfrac {26669714}{7875}b_{{101}}-\tfrac {433328537}{44625}b_{{110}}-\tfrac {225465802}{14875}b_{{200}},\nonumber
\end{eqnarray}
where $L_k^j$ denotes the homogeneous part of degree $j$ of $L_k$. We now turn to Theorem \ref{TeoBifHighOrder}. After a suitable change of variables in the perturbation parameters $(a_{002}, a_{011}, a_{020}, a_{101}, a_{110}, a_{200},$ $ b_{002}, b_{011}, b_{020})$, we can write 
$$L_i^1=u_i,\quad i=1,\dots,9.$$
Considering the new variables $\Lambda=(u_1,\dots,u_9,b_{101}, b_{110}, b_{200}, c_{002})$, we find a solution $\Lambda=s$ for the equations $L_{10}^2(\Lambda)=L_{11}^2(\Lambda)=L_{12}^2(\Lambda)=0$ such that $L_{13}^2(s)\neq 0$. Since
\begin{eqnarray}
L_{10}^2&=&-{\tfrac {
		42463\dots00071
		}{
		62328\dots85264
}}b_{{101}}b_{{200}}-{\tfrac {
		40320\dots92721
		}{
		31164\dots42632
}}b_{{110}}b_{{101}}-{\tfrac {
		65360\dots49715
		}{
		93493\dots27896
}}b_{{101}}c_{{002}}\nonumber\\
&&-{\tfrac {
		49690\dots18261
		}{
		12119\dots05468
}}b_{{110}}b_{{200}}-{\tfrac {
		91271\dots35547
		}{
		81806\dots36909
}}b_{{101}}^{2}-{\tfrac {
		46335\dots13063
		}{
		36358\dots16404
}}b_{{110}}^{2}\nonumber\\
&&+{\tfrac {
		58144\dots10977
		}{
		14543\dots65616
}}b_{{200}}^{2}-{\tfrac {
		27934\dots30905
		}{
		34627\dots15848
}}c_{{002}}b_{{200}}-{\tfrac {
		20443\dots06675
		}{
		13754\dots68284
}}c_{{002}}b_{{110}},\nonumber\\
L_{11}^2&=&{\tfrac {
		20939\dots26951
		}{
		30054\dots46160
}}b_{{101}}b_{{200}}+{\tfrac {
		19880\dots86461
		}{
		15027\dots23080
}}b_{{110}}b_{{101}}+{\tfrac {
		64454\dots32783
		}{
		90164\dots33848
}}b_{{101}}c_{{002}}\nonumber\\
&&+{\tfrac {
		40002\dots58387
		}{
		95411\dots04640
}}b_{{110}}b_{{200}}+{\tfrac {
		25716\dots92939
		}{
		22541\dots84620
}}{b_{{101}}}^{2}+{\tfrac {
		65259\dots07273
		}{
		50091\dots74360
}}{b_{{110}}}^{2}\nonumber\\
&&-{\tfrac {
		16384\dots09829
		}{
		40072\dots94880
}}{b_{{200}}}^{2}+{\tfrac {
		33058\dots21677
		}{
		40072\dots59488
}}c_{{002}}b_{{200}}+{\tfrac {
		57083\dots74365
		}{
		37899\dots60712
}}c_{{002}}b_{{110}},\nonumber\\
L_{12}^2&=&-{\tfrac {
		81825\dots04719
		}{
		11696\dots02500
}}b_{{101}}b_{{200}}-{\tfrac {
		14372\dots72644
		}{
		10818\dots73125
}}b_{{110}}b_{{101}}-{\tfrac {
		93191\dots38049
		}{
		12982\dots67750
}}b_{{101}}c_{{002}}\nonumber\\
&&-{\tfrac {
		16194\dots76093
		}{
		38467\dots60000
}}b_{{110}}b_{{200}}-{\tfrac {
		40196\dots79299
		}{
		35088\dots07500
}}{b_{{101}}}^{2}-{\tfrac {
		37739\dots55031
		}{
		28850\dots95000
}}{b_{{110}}}^{2}\nonumber\\
&&+{\tfrac {
		94762\dots70553
		}{
		23080\dots60000
}}{b_{{200}}}^{2}-{\tfrac {
		19118\dots35489
		}{
		23080\dots76000
}}c_{{002}}b_{{200}}-{\tfrac {
		65905\dots19663
		}{
		43657\dots69800
}}c_{{002}}b_{{110}},\nonumber
\end{eqnarray}
a solution $\Lambda=s$ is given by:
\begin{eqnarray}
b_{101}&=&1,\nonumber\\
b_{200}&=&\frac{\alpha}{3},\nonumber\\
b_{110}&=&{\tfrac{99478\dots56893{\alpha}^{2}-
		15595\dots75578\alpha-
		61109\dots90912}{
		74402\dots67872{\alpha}^{2}+
		13857\dots14244\alpha-
		17563\dots30264}},\nonumber\\
c_{002}&=&-{\tfrac {14326\dots21024{\alpha}^{4}+
		41672\dots01633{\alpha}^{3}-
		54197\dots40400{\alpha}^{2}-
		12687\dots07756\alpha+
		19077\dots08656}{
		65412\dots11040{\alpha}^{3}+
		29236\dots40360{\alpha}^{2}+
		16155\dots29280\alpha-
		40172\dots14560}},\nonumber
\end{eqnarray}
and
$$L_{13}^2(s)=\tfrac {60274\dots32675{\alpha}^{2}-
		33282\dots47954\alpha-
		87785\dots71976}{19996\dots30400},$$
where $\alpha$ is a real root of the cubic polynomial
$$37499\dots25384x^3 + 97562\dots80503x^2 - 13970\dots27210x - 27339\dots94856.$$

Note that $L_{13}^2(s)\neq 0$. Moreover, since the jacobian determinant

$$\det\dfrac{\partial (L_{10}^2(s),L_{11}^2(s),L_{12}^2(s))}{\partial (b_{110}, b_{200}, c_{002})}={\tfrac {15316\dots73503{\alpha}^{2}-
		11623\dots71962\alpha-
		78417\dots92584 }{
		20165\dots60960}}
,$$
is non-zero, we have that the hypersurfaces $L_{i}^2=0, i=10,11,12$ intersect transversally along $s$. Therefore, by Theorem \ref{TeoBifHighOrder}, there are perturbations of the center at the origin which yield $13$ limit cycles. \qed

\section{Acknowledgments}

We would like to thank Luiz Fernando Gouveia for the helpful insights and discussions which surely made this work richer. The first author is supported by S\~ao Paulo Research Foundation (FAPESP) grants 18/19726-5 and 19/10269-3. The second author is supported by S\~ao Paulo Research Foundation (FAPESP) grant 19/13040-7. The third author is partially supported by the Instituto Federal de Educa\c{c}\~ao, Ci\^encia e Tecnologia do Sul de Minas Gerais - IFSULDEMINAS.

\addcontentsline{toc}{chapter}{Bibliografia}
\bibliographystyle{siam}
\bibliography{Referencias.bib}
\end{document}